\documentclass[12pt,leqno]{article}
\usepackage{a4wide}
\usepackage{amsmath}
\usepackage{amsfonts}
\usepackage{amssymb}
\usepackage{theorem}
\usepackage{boxedminipage}
\usepackage{mathrsfs}
\usepackage{GV}
\usepackage{hyperref}
\usepackage[all]{xy}
\numberwithin{equation}{subsubsection}
\hypersetup{colorlinks=false}
\renewcommand{\modif}[1]{#1}
\date{}
\title{An extension of Greenberg's theorem to general valuation rings%
}
\author{Laurent Moret-Bailly
\thanks{The author is a member of the ANR project ``Points entiers et points rationnels''.
}
\medskip
\\
{\small IRMAR (Institut de Recherche Math\'{e}matique de Rennes, UMR 6625, CNRS)}\\
{\small Universit\'{e} de Rennes 1, campus de Beaulieu, F-35042 Rennes Cedex}\\
{\small \href{mailto:laurent.moret-bailly@univ-rennes1.fr}{laurent.moret-bailly@univ-rennes1.fr}}\\
{\small{\href{http://perso.univ-rennes1.fr/laurent.moret-bailly/}{http://perso.univ-rennes1.fr/laurent.moret-bailly/}}}
}
\begin{document}
%\selectlanguage{english}
\maketitle
\noindent\begin{boxedminipage}{\textwidth}
{\small{\slshape Cite as:}
L.\ Moret-Bailly, Manuscripta mathematica, 2011, DOI: 10.1007/s00229-011-0510-5.\\
The final publication is available at \href{http://www.springerlink.com}{www.springerlink.com}.}
\end{boxedminipage}
\begin{abstract}
We extend Greenberg's strong approximation theorem to schemes of finite presentation over valuation rings with arbitrary value group. As an application, we prove a closed image theorem (in the strong topology on rational points) for proper morphisms of varieties over valued fields.\smallskip\\
{\slshape 2010 Mathematics Subject Classification:}  13B40, %Étale and flat extensions; Henselization; Artin approximation
13F30, %valuation rings
12L10, %Ultraproducts
14G27. %Other nonalg. closed ground fields
\end{abstract}
%%%%%%%%%%%%%%%%%%%%%%%%%%%%%%%%%%
%\noindent\begin{boxedminipage}{\textwidth}
%{\small{\slshape Cite as:}
%L.\ Moret-Bailly, Manuscripta mathematica, 2011, DOI: 10.1007/s00229-011-0510-5.\\
%The final publication is available at \href{http://www.springerlink.com}{www.springerlink.com}.}
%\end{boxedminipage}
\smallskip
\begin{center}
\modif{\textsl{To Jan Denef, on the occasion of his 60th birthday}}
\end{center}%\vskip1cm
%\hfill\textsl{To Jan Denef, on the occasion of his 60th birthday}
\section{Introduction}
\Subsection{Notations}\label{notations}
Throughout this paper, we denote by $R$ a valuation ring,  by $K$ its fraction field, and by $\Gamma$ the valuation group (\modif{written} additively). The valuation is denoted by $\ord:K\to\Gi$. We put $\Gamma^+:=\{\alpha\in\Gamma\mid\alpha\geq0\}$.

The completion of $R$ is denoted by $\wh{R}$, with fraction field  $\wh{K}$; recall  that  $\wh{R}$ is a valuation ring with group $\Gamma$.

For each $\alpha\in\Gamma^+$, we put $I_{\alpha}:=\{x\in K\mid \ord\,{(x)}\geq\alpha\}$. This is a principal ideal of $R$, with quotient $R_{\alpha}:=R/I_{\alpha}$.

If  $X$ is an $R$-scheme, the sets $X(R_{\alpha})$ ($\alpha\in\Gamma^+$) form an inverse system, whose limit $\varprojlim_{\alpha\in\Gamma^+}X(R_{\alpha})$ is easily seen (although we shall not really need it)  to be $X(\wh{R})$. Indeed, this is immediate if $X$ is affine; in general, since each $R_{\alpha}$ is local, every element of the projective limit belongs to $\varprojlim_{\alpha\in\Gamma^+}V(R_{\alpha})$ where $V\subset X$ is an affine open subscheme.

Clearly, we have a natural map $X(R)\to\varprojlim_{\alpha\in\Gamma^+}X(R_{\alpha})$. \medskip

Our 
%first 
\modif{main result is the following:
\begin{thm}\textup{(strong approximation) }\label{main01} With the above notations, assume further  that $R$ is Henselian and that $\wh{K}$ is a separable extension of $K$. Let $X$ be an $R$-scheme of finite presentation. 

Then there exist a positive integer $N$ and an element $\delta\in\Gamma^+$ with the following property: for each $\alpha\in\Gamma^+$ and each $x\in X(R_{N\alpha+\delta})$, there is an $x'\in X(R)$ such that $x$ and $x'$ have the same image in $X(R_{\alpha})$.\smallskip\\
Equivalently:
\begin{equation}\label{EqTAF}
\forall\alpha\in\Gamma^+,\quad\im\bigl(X(R_{N\alpha+\delta})\to X(R_{\alpha})\bigr)=\im \bigl(X(R)\to X(R_{\alpha}) \bigr).
\end{equation}
\end{thm}
%
%This immediately implies:
\begin{subcor}\textup{(weak approximation) }\label{WA}
We keep the notations and assumptions of Theorem \rref{main01}. Then:
\begin{romlist}
\item\label{WA1} $X(R)$ is dense in $X(\wh{R})$ for the valuation topology.
\item\label{WA2} If $V$ is a $K$-scheme locally of finite type, then $V(K)$ is dense in $V(\wh{K})$ for the valuation topology.
\end{romlist}
\end{subcor}
\dem Theorem \ref{main01} immediately implies that $X(R)$ and $X(\wh{R})$ have the same image in $X(R_{\alpha})$, for each $\alpha\in\Gamma^+$. This proves \ref{WA1}. To prove \ref{WA2}, observe that $V(K)$ has an open covering by subsets of the form $j(\cU(R))$ where each $\cU$ is an affine $R$-scheme of finite presentation with an open immersion $j:\cU_{K}\inj V$ (see \ref{TopFacts}). Thus,  \ref{WA2} follows from  \ref{WA1} since $V(\wh{K})$ is then covered by the corresponding sets $\cU(\wh{R})$. \qed
\begin{subcor}\textup{(``infinitesimal Hasse principle'') }\label{PHI}
With the notations and assumptions of Theorem \rref{main01}, we have the equivalence:
$$X(R)\neq\emptyset\quad \Iff\quad \forall\gamma\in\Gamma^+, \:X(R_{\gamma})\neq\emptyset.$$
\end{subcor}
\dem Taking $\alpha=0$ in (\ref{EqTAF}), we see that if $X(R_{\delta})\neq\emptyset$ then $X(R)\neq\emptyset$ .\qed
}\bigskip

The author's original motivation for proving these results is that they have deep consequences for the topology of varieties over valued fields. For instance, as an easy consequence of Corollary \ref{PHI}, we obtain:

\begin{thm}\label{ThPropre} Assume $R$ is Henselian and $\wh{K}$ is separable over $K$. Let $f:X\to Y$ be a \emph{proper} morphism of $K$-schemes of finite type. Then the induced map $f_{K}:X(K)\to Y(K)$ has closed image (for the topology defined by the valuation).
\end{thm}
\begin{subrem}
Theorem \ref{ThPropre} is of course trivial if $K$ is a \emph{local} field (i.e.\ locally compact), since $f_{K}$ is then a proper map. But apart from this case, and even if $R$ is a discrete valuation ring, $f_{K}$ is not a closed map in general.
\end{subrem}

\Subsection{Related results}
Theorem \ref{main01} of course generalizes Greenberg's strong approximation theorem \cite{Green66}, which is the special case where $R$ is a discrete valuation ring (the separability of $\wh{K}$ meaning in this case that $R$ is excellent). In fact, Greenberg's original proof extends rather easily to valuation rings of height one provided the fraction field has characteristic zero.

The method used here is due to Becker, Denef, Lipshitz and van den Dries \cite{BeDeLivdD79}: in fact, most of our proof is shamelessly copied from there,  with the exception of the separability property \ref{regular}\,\ref{regular2} which is proved in \cite{BeDeLivdD79} by a ramification index argument which breaks down for nondiscrete valuations. 

Similar methods are used in \cite{BeDeLivdD79}, and also by Denef and Lipshitz in \cite{DeLi80} to obtain strong approximation theorems  more general than Greenberg's (of the kind considered by Artin, Popescu and others); typically, these are derived from the corresponding ``weak'' approximation theorems. The ground rings in these results are subrings of power series rings over discrete valuation rings.

Schoutens \cite[Theorem 2.4.1]{Schou88} has proved the weak approximation theorem \ref{WA} for certain subrings of $A[[T_{1},\dots, T_{n}]]$ where $A$ is a complete valuation ring of height one.

The approach of Elkik \cite{Elkik73} leads to strong approximation results without excellence assumption on the base ring, but with (generic) smoothness assumptions on the scheme $X$. The extension to certain non-Noetherian bases (including Henselian valuation rings of height one), outlined in \cite[Remarque 2, p.~587]{Elkik73}, is carried out   in \cite[1.16]{Abbes11}.

\Subsection{Organization of the paper}
In Section \ref{SecBasic} we review some basic facts about ultraproducts, in particular about ultrapowers of $R$. We then explain how Theorem \ref{main01} reduces to two technical results involving such ultrapowers (namely, the separability theorem \ref{regular}\,\ref{regular2} and the lifting theorem \ref{mainbis}). This reduction is the ``formal'' part of the proof. Theorems \ref{regular} and \ref{mainbis} are proved in section  \ref{sec:proof}, and Theorem \ref{ThPropre} in section \ref{sec:dempropre}. 

\Subsection{Acknowledgments}

The author is grateful to Hans Schoutens for pointing out his work \cite{Schou88}, \modif{to Ofer Gabber who removed an annoying assumption from the main result}, and  to Jan Denef for his strategic help. \modif{He also thanks the referee for his remarks.}

%%%%%%%%%%%%%%%%%%%%%%%%%%%%%%%%%%%%
\section{Basic constructions}\label{SecBasic}
We keep the notations of the introduction.
\Subsection{Ultraproducts: basic definitions}
(For details on these constructions, see \cite[\S2]{BeDeLivdD79} or \cite{Schou2010}). 

%\Subsubsection{Ultraproducts. }
We fix an infinite set $\W$ and a nonprincipal  ultrafilter $\cU$  of subsets of $\W$. We shall  say that a property $P(w)$ holds ``for almost all $w\in W$''  if the corresponding subset of $\W$ belongs to $\cU$.

If $A=(A_{w})_{w\in \W}$ is a family of algebraic structures (sets, rings, groups, ordered groups\dots) indexed by $\W$, we denote by $A_{\ult}$ (notation taken from \cite{Schou2010}) the corresponding ultraproduct. It is \emph{usually} defined as  the quotient of $\prod_{w \in \W}A_{w}$ by the equivalence relation  ``equality almost everywhere''. It will often be convenient to use the more explicit notation (inspired from the same source)  $\ulim_{\cU,\,w}A_{w}$ for $A_{\ult}$.

However, the above definition works as expected only if either the sets $A_{w}$ are nonempty, or almost all are empty. In general, the correct definition (which we shall use here)  is
\begin{equation}\label{EqDefUltraprod}
\ulim_{\cU,\,w}A_{w}:=\varinjlim\limits_{U\in\cU}\:\prod_{w\in U}A_{w}
\end{equation}
where $\cU$ is ordered by reverse inclusion and the transition maps are the obvious projections.

Given $U\in\cU$ and an element $x=(x_{w})_{w \in U}\in \prod_{w \in U}A_{w}$, we denote its class in $A_\ult$ by
% $x_{\ult}$, or 
$[x_{w}]_{w \in U}$, or 
$\ulim_{\cU,\,w}x_{w}$.  

If the family is constant ($A_{w}=A$, independent of $w$) we obtain the $\cU$-\emph{ultrapower} of $A$, denoted by
\begin{equation}\label{EqDefUltrapuiss}
\upw_{\cU}A=\ulim_{\cU,w}A.
\end{equation}

\Subsection{Ultraproducts and the functor of points}
If $A$ is a ring and $Y$ is an $A$-scheme, we wish to know whether the functor of points $B\mapsto Y(B)$ (from $A$-algebras to sets) ``commutes with ultraproducts''. Given the definition of ultraproducts, this must involve compatibility properties of this functor  with products and with direct limits. Now, recall the following facts:

\begin{subprop}\textup{(special case of {\cite[(8.8.2)]{EGA4_III}}) }\label{PointsLimInd} Let $A$ be a ring,  $(B_{\lambda})_{\lambda\in\Lambda}$  a filtering inductive system of $A$-algebras, and $Y$ an $A$-scheme. Consider the natural map
\begin{equation}\label{EqPointsLimInd}
\alpha:\varinjlim_{\lambda\in\Lambda}Y(B_{\lambda})\ffl Y\bigl(\varinjlim_{\lambda\in\Lambda}B_{\lambda}\bigr).
\end{equation}
If $Y$ is locally of finite type (resp. locally of finite presentation) over $A$, then $\alpha$  is injective (resp. bijective).\qed
\end{subprop}

\begin{subprop}\label{PointsProd} Let $A$ be a ring,  $(B_{i})_{i\in I}$  a family of $A$-algebras, and $Y$ an $A$-scheme. Consider the natural map
\begin{equation}\label{EqPointsProd}
\beta:Y\Bigl(\prod_{i\in I}B_{i}\Bigr)\ffl \prod_{i\in I}Y(B_{i}).
\end{equation}
\begin{romlist}
\item\label{PointsProd1} If $Y$ is affine, $\beta$ is bijective.
\item\label{PointsProd2} If $Y$ is quasiseparated, $\beta$ is injective.
\item\label{PointsProd3} If $Y$ is quasicompact and quasiseparated, \emph{and} each $B_{i}$ is a local ring, then $\beta$ is bijective.
\end{romlist}
\end{subprop}
\dem \ref{PointsProd1} is immediate since the $\Spec$ functor takes arbitrary products of rings to sums in the category of \emph{affine} schemes. 

A proof of \ref{PointsProd2} and \ref{PointsProd3} can be found embedded in the proof of \cite[Theorem 3.6]{ConAdelic}. Note that \ref{PointsProd3} appears in \cite[Lemme 3.2]{Oe84}, although the quasiseparated assumption is missing there.\qed
\smallskip

Now, let  $A$ be a ring, and let $(A_{w})_{w\in W}$ be a family of $A$-algebras indexed by $W$, with ultraproduct $A_{\ult}$. If $Y$ is an $A$-scheme, we want to compare the sets $Y(A_\ult)$ and $\ulim_{\cU,w}Y(A_{w})$. 

For each $U\in\cU$, put $A_{U}:=\prod_{w\in U}A_{w}$. Since $A_{\ult}=\varinjlim_{U}A_{U}$, we have a natural map $\alpha:\varinjlim_{U\in\cU}Y(A_{U})\to Y(A_{\ult})$, to which Proposition \ref{PointsLimInd} applies.

On the other hand, for each $U$ we have a map $\beta_{U}: Y(A_{U})\to \prod_{w\in U}Y(A_{w})$ (of the type considered in Proposition \ref{PointsProd}) and, passing to the limit, a map $\beta: \varinjlim_{U\in\cU}Y(A_{U})\to \ulim_{\cU,w}Y(A_{w})$. Finally, we have constructed a diagram of sets
\begin{equation}\label{DiagCommUlt}
\parbox{10cm}{$\xymatrix{&\varinjlim\limits_{U\in\cU}Y(A_{U})\ar@<1ex>[r]^(0.4){\alpha}\ar[d]_{\beta}
& Y(\varinjlim\limits_{U\in\cU}A_{U})=Y(A_{\ult})\\
\ulim_{\cU,w}Y(A_{w})\ar@<.6ex>@{=}[r]&\varinjlim\limits_{U\in\cU}\prod\limits_{{w\in U}}Y(A_{w}).
}$}
\end{equation}
Combining Propositions \ref{PointsLimInd} and \ref{PointsProd}, we obtain:
\begin{subprop}\label{PointsUltraprod} With the above notations and assumptions, assume that: 
\begin{sitemize}
\item  $Y$ is finitely presented over $A$, and
\item $Y$ is affine, or each $A_{w}$ is a local ring.
\end{sitemize}
Then the maps $\alpha$ and $\beta$ in \textup{(\ref{DiagCommUlt})} are bijective. In particular, we have a natural bijection
$$\ulim_{\cU,w}Y(A_{w})\fflis Y(A_{\ult}).\eqno{\qed}$$
\end{subprop}

\begin{subrems}{\ }
\begin{remlist}
\item In the affine case, $Y(B)$ (for an $A$-algebra $B$) is the set of $B$-valued solutions of a given finite system of polynomial equations with coefficients in $A$. Thus, Proposition \ref{PointsUltraprod} in this case will be seen by model theorists as an instance of \L{o}\'{s}' theorem.
\item There is no ``direct'' map between the two sets in \ref{PointsUltraprod}, valid for \emph{all} $A$-schemes $Y$. To see this, take for $W$ the set of prime numbers (and for $\cU$ any nonprincipal ultrafilter), and take $A=\ZZ$, $A_{p}=\FF_{p}$ for all $p\in W$. In particular, $A_{\ult}$ is a field of characteristic zero.

First, try  $Y=\Spec\,\QQ$. Then $Y(A_{\ult})$ has one element, while $Y(A_{p})=\emptyset$ for each $p\in W$. Hence in this case there is no map \modif{from $Y(A_{\ult})$ to $\ulim\nolimits_{\cU,p}Y(A_{p})$}. 

Now, take $Y=\coprod_{p\in W}\Spec\,\FF_{p}$. Then $Y(A_{\ult})=\emptyset$, and $Y(A_{p})$ has one element for each $p$, so $\ulim \nolimits_{\cU,p}Y(A_{p})$ also has one element and there is no map \modif{from $\ulim\nolimits_{\cU,p}Y(A_{p})$ to  $Y(A_{\ult})$}. 
\item Proposition \ref{PointsUltraprod} would fail with the traditional definition of ultraproducts: it may happen that $Y(A_{\ult})\neq\emptyset$ but $Y(A_w)=\emptyset$ for some $w\in W$.
\item With the assumptions of \ref{PointsUltraprod}, we have in particular:
\begin{equation}\label{eq:PolEq2}
Y\Bigl(\ulim_{\cU,\,w}A_{w}\Bigr)\neq\emptyset\quad\Leftrightarrow\quad Y(A_{w})\neq\emptyset\text{ for almost all $w$}.
\end{equation}
An interesting special case is when $A_{w}=B$, a fixed local $A$-algebra: we then have the equivalence
\begin{equation}\label{eq:PolEq3}
Y\Bigl(\upw_{\cU}(B)\Bigr)\neq\emptyset\quad\Leftrightarrow\quad Y(B)\neq\emptyset.
\end{equation}
\end{remlist}
\end{subrems}

\Subsection{Ultrapowers of valuation rings}
Take our valuation ring $R$, and consider $R_{\ult}=\upw_{\cU}R$. \modif{Using \L{o}\'{s}' theorem, one checks that} this is a valuation ring with fraction field $K_\ult=\upw_{\cU}K$ and valuation group $\Gamma_\ult=\upw_{\cU}\Gamma$; moreover, $R_\ult$ is Henselian if $R$ is.

 We have canonical embeddings $R\inj R_\ult$, $\Gamma\inj\Gamma_\ult$. We shall denote the valuation on $K_\ult$ by $\ord_\ult:K_\ult\to\Gamma_\ult\cup\{\infty\}$. Thus, if $z=(z_{w})_{w \in \W}\in K^\W$, we have $\ord_\ult [z_{w}]_{w \in \W}=[\ord(z_{w})]_{w \in \W}$.
\Subsubsection{Some quotients of $R_{\ult}$: principal ideals. }\label{PpalIdeals}
Each element ${\alpha_{\ult}}$ of $\Gamma^+_{\ult}$ defines  a principal ideal $I_{\alpha_{\ult}}\subset R_{\ult}$ and a quotient ring $(R_{\ult})_{\alpha_{\ult}}$, which we denote by $R_{\ult,\alpha_{\ult}}$. If we write $\alpha_{\ult}=[\alpha_{w}]_{w\in W}$ for some family $(\alpha_{w})\in (\Gamma^+)^W$, we immediately check that the canonical surjection $R^W\surj \prod_{w\in W}R_{\alpha_{w}}$ induces an isomorphism 
$$ R_{\ult,\alpha_{\ult}} \fflis \ulim_{\cU,\,w}R_{\alpha_{w}}.$$

\Subsubsection{Some quotients of $R_{\ult}$: prime ideals. }\label{sub:primes}
Recall that a subset $C$ of an ordered set $(S,\leq)$ is \emph{convex} if whenever $a\in C$, $b\in C$, $x\in S$ and $a\leq x\leq b$, then $x\in C$. (Convex subgroups of a totally ordered group are called \emph{isolated}  in \cite{BourAC5-6}).

Let $C$ be a convex subgroup of $\Gamma_\ult $. We denote by $P_{C}\subset R_\ult $ the ideal
$$%\begin{array}{rcl}
P_{C}\;=\;\left\{x\in R_{\ult} \mid \ord_{\ult}(x)\notin C\right\}\;=\;\left\{x\in R_{\ult} \mid \ord_{\ult}(x)> C\right\},
%\end{array}
$$
the latter condition meaning of course that $\ord_{\ult}(x)> \alpha$ for all $\alpha\in C$.
This is  a prime ideal of $R_\ult $ (they are all of this form), and the quotient $R^\ang{C}:=R_\ult /P_{C}$ is a valuation ring with group $C$. If $C$ contains $\Gamma$, the canonical map $R\to R_\ult \to R^\ang{C} $ is injective. 

The ideal $P_{C}$ is not principal in general, but it is  the (totally ordered) union of the principal ideals contained in it. So we can write (as $R_{\ult}$-algebras)
$$R^\ang{C}=\varinjlim_{\alpha_{\ult}>C}R_{\ult,\alpha_{\ult}}.$$
We shall be interested only in convex subgroups $C$ satisfying $\Gamma\subset C \underset{\neq}{\subset} \Gamma_{\ult}$. (If $C=\Gamma_{\ult}$, then $P_{C}=\{0\}$, and the above direct limit runs over the empty set; otherwise we have an honest filtering colimit). Any such subgroup contains the convex hull $\Gamma_{c}$ of $\Gamma$ in $\Gamma_\ult$. We can think of $P_{\Gamma_{c}}$ as the ideal of elements of $R_\ult $ with ``infinitely large'' valuation. 

Thus we have a diagram of valuation rings
$$
\xymatrix{R_\ult \ar@{>>}[r]&R_{\ }^\ang{\Gamma_{c}}\\
R^{\strut}\ar@{^{(}->}[u]\ar@{^{(}->}[ur]\\
}$$
and our general strategy for solving equations over $R$ will be ``find solutions in $R^\ang{\Gamma_{c}}$, lift them to $R_\ult $, and then extract solutions in $R$''. As we shall see now, the first and third steps are essentially trivial.

\begin{subprop}\label{Los}
Let $C$ be a \emph{proper} convex subgroup of $\Gamma_{\ult}$ containing $\Gamma$, and let $X$ be an $R$-scheme of finite presentation. Then we have the implications
$$
X(R)\neq\emptyset \quad\Iff \quad X(R_\ult )\neq\emptyset \quad\Impl\quad X(R^{\ang{C}})\neq\emptyset \quad \Iff\quad  \forall\alpha\in\Gamma^+,X(R_{\alpha})\neq\emptyset.
$$ 
\end{subprop}
\pf The first two ``$\Impl$'' are obvious, and the first equivalence follows from  \eqref{eq:PolEq3}.

Assume \modif{$X(R^{\ang{C}})\neq\emptyset$}, and take $\alpha\in\Gamma^+$. The ideal \modif{$P_{C}$} is  contained \modif{in $P_{\Gamma_{c}}$, hence in}  $I_{\alpha}R_{\ult}$, and therefore $R_{\ult}/I_{\alpha}R_{\ult}$ is a quotient of \modif{$R^{\ang{C}}$}, which implies that $X(R_{\ult}/I_{\alpha}R_{\ult})\neq\emptyset$. But $R_{\ult}/I_{\alpha}R_{\ult}$ is immediately seen to be the ultrapower $\upw_{\cU}R_{\alpha}$, whence $X(R_{\alpha})\neq\emptyset$ by \eqref{eq:PolEq3} again. This proves the last ``$\Impl$''.

%(\ref{Los2}) 
Finally, assume the last condition in the chain. Since $C \underset{\neq}{\subset} \Gamma_{\ult}$ by assumption, we can pick some $\alpha_{\ult}>C$ in $\Gamma_{\ult}$, and it suffices to show that $X(R_{\ult,\alpha_{\ult}})\neq\emptyset$ since $R^\ang{C}$ is a quotient of $R_{\ult,\alpha_{\ult}}$. Now represent $\alpha_{\ult}$ as $\ulim\limits_{\cU,w}\alpha_{w}$ for some $(\alpha_{w})\in \Gamma^W$: then from \ref{PpalIdeals} we have $R_{\ult,\alpha_{\ult}}=\ulim\limits_{\cU,w}R_{\alpha_{w}}$, whence, by \ref{PointsUltraprod}, $X(R_{\ult,\alpha_{\ult}})=\ulim\limits_{\cU,w}X(R_{\alpha_{w}})\neq\emptyset$ since each $X(R_{\alpha_{w}})$ is nonempty.\qed

\begin{subrem} It may happen that $\Gamma_{c}=\Gamma_{\ult}$, in which case there is no $C$ as in the proposition. This is the case in particular if $W$ is ``too small'' in the sense that $\Gamma$ has no cofinal subset of cardinality $\leq\mathrm{Card\,}W$.

On the other hand, if we restrict ourselves to those $(W,\cU)$ such that $\Gamma_{c}\neq\Gamma_{\ult}$, then \modif{\ref{Los} shows that the condition $X(R^{\ang{C}})\neq\emptyset$  is equivalent to $X(R^{\ang{\Gamma_{c}}})\neq\emptyset$, hence independent of $C$ (and  even independent of $(W,\cU)$, subject to the above restriction).}
%In the last part of the proof, we have in fact constructed a natural map from $\prod_{w \in \W}X(R_{w})$ to $X(R^\ang{\Gamma_{c}})$.
\end{subrem}\medskip

Let us now state the technical results from which  \ref{main01} will be derived. First, a structure theorem for the fraction fields of the rings $R^\ang{C}$:

\begin{thm}\label{regular} Let $C$ be a convex subgroup of $\Gamma_{\ult}$ containing $\Gamma$. Consider the extension $K^\ang{C}:=\Frac(R^\ang{C})$ of $K$.
\begin{romlist}
\item \label{regular1}  If $R$ is \emph{complete}, then  $K^\ang{C}$ is a \emph{regular} extension of $K$, i.e.\ $K^\ang{C}$ is linearly disjoint from every finite extension of $K$. (In other words, $K^\ang{C}$ is a geometrically integral $K$-algebra).
\item \label{regular2}  If $\wh{K}$ is separable over $K$, then so is  $K^\ang{C}$.  (In other words, $K^\ang{C}$ is a geometrically reduced $K$-algebra).
\item  \label{regular3}  \modif{If $K$ is separably closed in $\wh{K}$ \emph{(e.g.\ if $R$ is Henselian \cite[{F, Th.\ 4, Cor.\ 2 p.\ 190}]{Rib68})}, then it is separably closed in   $K^\ang{C}$. (In other words, $K^\ang{C}$  is a \emph{primary} extension of $K$, or equivalently a geometrically connected $K$-algebra).}
\end{romlist}
\end{thm}

\modif{A  word of warning may be appropriate here:  for a valuation ring, ``complete'' does \emph{not} imply ``Henselian'', except if the value group $\Gamma$ has height one, i.e.\ is isomorphic to a subgroup of $\RR$ with the induced ordering. }\smallskip

Theorem \ref{regular} will be proved in section \ref{sec:proof}. As we shall see, assertions \ref{regular2} and \ref{regular3} are easy consequences of \ref{regular1}. For us, the useful one is the separability property \ref{regular2}, which will be used, also in section \ref{sec:proof}, to prove the following result:

\begin{thm}\textup{(Lifting theorem) }\label{mainbis} Assume that $R$ is Henselian and that  $\wh{K}$ is separable over $K$. For each convex subgroup $C\subset\Gamma_{\ult}$ containing $\Gamma$, the canonical map $X(R_\ult )\to X\left(R^\ang{C}\right)$ is onto.
\end{thm}

We shall end this section by deducing Theorem \ref{main01} from the lifting theorem.

\modif{\Subsection{Proof of Theorem \ref{main01} (from Theorem \ref{mainbis})}
We argue by contradiction. Thus, assume that for all $N\in\ZZ_{>0}$ and $\delta\in\Gamma^+$ there exist $\alpha_{N,\delta}\in\Gamma^+$ and $\xi_{N,\delta}\in X(R_{N\alpha_{N,\delta}+\delta})$ such that the image of $\xi_{N,\delta}$ in $X(R_{\alpha})$ does not lift to $X(R)$. Using the axiom of choice we fix such families $\left(\alpha_{N,\delta}\right)$ and $\left(\xi_{N,\delta}\right)$.
%; we may assume $\alpha_{N,\delta}>0$ for each $\alpha$ by increasing them if necessary (this is possible unless $X(R_{\alpha})=\emptyset$ for some $\alpha$, in which case the result is trivial).
For simplicity, put $\beta_{N,\delta}:=N\alpha_{N,\delta}+\delta$. 

Now, let us choose our ultrafilter: we take $W:=\ZZ_{>0}\times \Gamma^+$, and pick an ultrafilter $\cU$ on $W$, containing all the sets $w+W$ ($w\in W$). The ultrapowers $\ZZ_{\ult}$ and $\Gamma_{\ult}$ contain in particular the ``diagonal'' elements
$$H:=\ulim_{\cU,(N,\delta)}N,\qquad \Delta:=\ulim_{\cU,(N,\delta)}\delta$$
and our choice of $\cU$ implies that 
\begin{equation}\label{EqDemMain01:2}
H>\ZZ \text{ (in $\ZZ_{\ult}$)\quad and \quad}\Delta>\Gamma\text{ (in $\Gamma_{\ult}$)}. 
\end{equation}
For $w=(N,\delta)\in W$, we of course write $\alpha_{w}$ for $\alpha_{N,\delta}$. Now we consider the elements
$$\begin{array}{rcll}
\alpha_{\ult}&:=&\ulim\limits_{\cU,w}\alpha_{w}&\in \Gamma_{\ult}^+,\\
\beta_{\ult}&:=&\ulim\limits_{\cU,w}\beta_{w}&\in \Gamma_{\ult}^+,\\
\xi_{\ult}&:=&\ulim\limits_{\cU,w}\xi_{w}&\in \ulim_{\cU,w}X\left(R_{\beta_{w}}\right)\cong X\Bigl(\ulim_{\cU,w}R_{\beta_{w}}\Bigr)=X\left(R_{\ult,\beta_{\ult}}\right)
\end{array}
$$
where in the last line we have used \ref{PointsUltraprod} and \ref{PpalIdeals}. An equivalent definition of $\beta_{\ult}$ is of course
$\beta_{\ult}=H\alpha_{\ult}+\Delta$
(note that  $\Gamma_{\ult}$ is an ordered $\ZZ_{\ult}$-module in a natural way); in particular, from (\ref{EqDemMain01:2}) (and since $\alpha_{\ult}\geq0 $) we see that $\beta_{\ult}\geq\ZZ\alpha_{\ult}+\Gamma$ in $\Gamma_{\ult}$, and consequently
$$\beta_{\ult}\geq C:=\text{ convex hull of $\ZZ\alpha_{\ult}+\Gamma$ in $\Gamma_{\ult}$}.$$
This means that $R^\ang{C}$ is a quotient of $R_{\ult,\beta_{\ult}}$. In turn, $R_{\ult, \alpha_{\ult}}$ is a quotient of $R^\ang{C}$, by definition of $C$. Thus we have a diagram of sets
$$\xymatrix{&X(R_{\ult})\ar[d]\\
\xi_{\ult}\in X(R_{\ult,\beta_{\ult}})\ar[r]&X(R^\ang{C})\ar[r]& X(R_{\ult, \alpha_{\ult}})
}$$
By Theorem \ref{mainbis}, the image of $\xi_{\ult}$ in $R^\ang{C}$ lifts to an element $\eta_{\ult}\in X(R_{\ult})$. By construction, $\xi_{\ult}$ and $\eta_{\ult}$ have the same image in $X(R_{\ult, \alpha_{\ult}})$. This means, using \ref{PpalIdeals}, that, for almost all $w\in W$, the image of $\xi_{w}$ in $X(R_{\alpha_{w}})$ lifts to $X(R)$. This contradicts our initial choices.\qed}

\section{Proof of Theorems \ref{regular} and  \ref{mainbis}}\label{sec:proof}
\Subsection{Finite extensions}\label{Ssec:finiteExt}
In this section we assume $R$ \emph{complete}. Let $K_{1}$ be a finite extension of $K$, of degree $d$. Then the valuation $\ord$ has an extension $\ord_{1}$ to $K_{1}$, with group $\Gamma_{1}\supset\Gamma$ (this is true for any extension). Moreover, we know that the index $(\Gamma_{1}:\Gamma)$ is finite (and in fact $\leq d$) \cite[VI, \S\,8, $\mathrm{n}^\circ$\, 3, th.~1]{BourAC5-6}; in particular, $\Gamma$ is cofinal in $\Gamma_{1}$. 

(Note that, unless it has height $1$,  $\ord$ may have several extensions to $K_{1}$; however, they are all dependent, i.e.~they define the same topology on $K_{1}$ \cite[\S\,8, $\mathrm{n}^\circ$\, 2, cor.~1]{BourAC5-6}).

We denote by $R_{1}\subset K_{1}$ the ring of $\ord_{1}$. The situation is complicated by the fact that $R_{1}$ is not necessarily a finitely generated  $R$-module. To address this, we shall define  substitutes for $R_{1}$ and  $\ord_{1}$ as follows: choose a $K$-basis $\cB$ of $K_{1}$ whose elements are integral over $K$ (hence $\cB\subset R_{1}$), and put $R_{0}:=R[\cB]\subset R_{1}$. Then $R_{0}$ is a finite $R$-algebra with fraction field $K_{1}$. Since $R$ is a valuation ring, $R_{0}$ is a \emph{free} $R$-module of rank $d$, so we can fix a basis $\cB_{0}=(e_{1}=1,e_{2},\dots,e_{d})$ of $R_{0}$ over $R$. Now, for each $z=\sum_{i=1}^d x_{i}e_{i} \in K_{1}$ (with the $x_{i}$'s in $K$) we can put 
$$\ord_{0}\,(z):=\min_{i}\ord\,(x_{i}).$$
\begin{sublem}\label{2topologies} With the above assumptions and notations, the function
$$f:=\ord_{1}-\ord_{0}:K_{1}^*\ffl\Gamma_{1}$$
is bounded.
%there is an element $c\in\Gamma$ withthe following property:
%$$\forall z \in K_{1},\quad  \ord_{0}(z)\leq\ord_{1}(z)\leq \ord_{0}(z)+c.$$
\end{sublem}
\pf 
Let us introduce the ``balls''
$$\begin{array}{rcll}
B_{0}(\alpha)&:=&\left\{z\in K_{1}\mid \ord_{0}(z)\geq\alpha\right\}&(\alpha\in\Gamma)\\
B_{1}(\alpha)&:=&\left\{z\in K_{1}\mid \ord_{1}(z)\geq\alpha\right\} &(\alpha\in\Gamma_{1}).
\end{array}
$$
Thus, we have $B_{0}(0)=R_{0}$ and $B_{1}(0)=R_{1}$. The family $\left(B_{1}(\alpha)\right)_{\alpha\in\Gamma_{1}}$ (resp. $\left(B_{0}(\alpha)\right)_{\alpha\in\Gamma}$) is a basis of neighbourhoods of $0$ in $K_{1}$ for the topology defined by $\ord_{1}$ (resp.~for the product topology on $K_{1}$, identified with $K^d$ via $\cB_{0}$); since $\Gamma$ is cofinal in $\Gamma_{1}$ we can even restrict the first family to $\Gamma$. Note that we have $t\,B_{0}(\alpha)=B_{0}(\alpha+\ord(t))$ for $\alpha\in\Gamma$ and $t\in K$, and similarly for $B_{1}$.

Since $K$ is complete, our two topologies are in fact the same \cite[chap.~6, \S5, $\mathrm{n}^\circ$ 2, prop. 4]{BourAC5-6}. In fact, we trivially have $B_{0}(\alpha)\subset B_{1}(\alpha)$ for all $\alpha$ (look at the definition of $\ord_{0}$); but by \cite{BourAC5-6} we also have, say, $B_{1}(\lambda)\subset B_{0}(0)$ for some $\lambda\in\Gamma$, whence $B_{1}(0)\subset B_{0}(-\lambda)$ and, by scaling, 
$$ B_{0}(\alpha)\subset B_{1}(\alpha)\subset B_{0}(\alpha-\lambda)$$
for all $\alpha\in\Gamma$. 

Returning to the function $f$, observe that $f(tz)=f(z)$ for $t\in K^*$. Next, since $\Gamma$ has finite index in $\Gamma_{1}$, there is a finite subset $\Sigma\subset\Gamma_{1}$ such that each $z\in K_{1}^*$ can be written $z=tz_{1}$ with $t\in K$ and $\ord_{1}(z_{1})\in \Sigma$. It follows that it is enough to bound $f(z)$ whenever $z$ is in the ``annulus'' $U:=B_{1}(r)\setminus B_{1}(r')$, for  any fixed  $r\leq\min({\Sigma})$ and $r'>\max(\Sigma)$, which we may (and do) take in $\Gamma$.  Now from the above inclusions we have $B_{0}(r')\subset B_{1}(r')\subset B_{1}(r)\subset B_{0}(r-\lambda)$, whence $U\subset B_{0}(r-\lambda)\setminus B_{0}(r')$. In other words, $\ord_{0}$ (hence also $f$) is bounded on $U$, which completes the proof.\qed

\Subsection{Proof of Theorem \ref{regular}}
We adopt the notations and assumptions of Theorem \ref{regular}. Let us first show that assertion \ref{regular1} easily implies  \ref{regular2} and  \ref{regular3}. Consider  the commutative diagram of \modif{fraction fields:
$$\begin{array}{ccccl}
K& \inj & K_\ult  & \inj & K^\ang{C}\\
\cap && \cap && \cap\\
\wh{K}& \inj & \wh{K}_\ult  & \inj & \wh{K}^\ang{C}.
\end{array}
$$}
Assuming \ref{regular1} (applied to $\wh{R}$), the extension of fraction fields $\wh{K}\inj\wh{K}^\ang{C}:=\Frac\,(\wh{R}^\ang{C})$ (bottom line) is regular. This proves that $\wh{K}^\ang{C}/K$ is separable (resp. primary) if $\wh{K}/K$ is, and the same holds for the subextension ${K}^\ang{C}/K$. \modif{This implies \ref{regular2} and  \ref{regular3}, as promised}.

From now on, we assume $R$ complete. Let $K_{1}$ be a finite extension of $K$: we need to prove that  $K^\ang{C}$ and $K_{1}$ are linearly disjoint over $K$. Put $d:=[K_{1}:K]$. 

We now apply the constructions (and keep the notations) of \ref{Ssec:finiteExt}. We also have ultrapowers $R_{0,\ult}\subset R_{1,\ult}\subset K_{1,\ult}$, and a valuation $\ord_{1,\ult}$ on  $K_{1,\ult}$ with ring  $R_{1,\ult}$ and group  $\Gamma_{1,\ult}$. Since $\Gamma\subset\Gamma_{1}$ has finite index, so does $\Gamma_{\ult}\subset\Gamma_{1,\ult}$. We denote by $C_{1}$ the convex hull of $C$ in $\Gamma_{1,\ult}$ (a convex subgroup containing $\Gamma_{1}$): this defines a prime ideal $P_{C_{1}}$ of $R_{1,\ult}$, with quotient $R_{1}^{\ang{C_{1}}}$. We immediately see that $C=\Gamma_{\ult}\cap C_{1}$ (since $C$ is already convex in $\Gamma_{\ult}$) and $P_{C}=R_{\ult}\cap P_{C_{1}}$.

The following lemma clearly implies that the fields $K_{1}=\Frac(R_{0})$ and $K^\ang{C}=\Frac(R^\ang{C})$ are linearly disjoint over $K$, thus completing the proof of \ref{regular}\,\ref{regular1}:
\begin{sublem}\label{domain} $R_{0}\otimes_{R}R^\ang{C}$ is an integral domain.
\end{sublem}
\pf Since $R_{0}$ is finite free over $R$, the natural map $R_{0}\otimes_{R}R^\W\to R_{0}^\W$ is an isomorphism, and one readily checks that the same holds for  $R_{0}\otimes_{R}R_\ult \to R_{0,\ult}$. Hence $R_{0}\otimes_{R}R^\ang{C}$ is isomorphic to $R_{0,\ult}/P_{C}R_{0,\ult}$. So, our task is to show that $P_{C}R_{0,\ult}$ is a prime ideal of $R_{0,\ult}$. In fact we shall prove that $P_{C}R_{0,\ult}=R_{0,\ult}\cap P_{C_{1}}$, which implies the claim since $P_{C_{1}}\subset  R_{1,\ult}$ is prime.

Since $(e_{1},\dots,e_{d})$ is a basis of $R_{0,\ult}$ over $R_\ult $, every element of $R_{0,\ult}$ can be written as 
$$x=\sum_{i=1}^d x^{(i)}_{\ult}\:e_{i} \quad (x^{(i)}_{\ult}\in R_{\ult})$$
and this $x$ is in $P_{C}R_{0,\ult}$ if and only if each coordinate $x^{(i)}_{\ult}$ is in $P_{C}$, or equivalently in $P_{C_{1}}$. In other words:
$$\begin{array}{rcl}
x_{\ult}\in P_{C}R_{0,\ult}\quad&\Iff&\forall i\in\{1,\dots,d\}, \ord_{\ult}\,(x^{(i)}_{\ult})>C_{1}\\
&\Iff& \min\limits_{1\leq i\leq d}\ord_{\ult}\,(x^{(i)}_{\ult})>C_{1}\\
&\Iff&\ord_{0,\ult}\,(x_{\ult})>C_{1}.
\end{array}$$
But it follows from Lemma \ref{2topologies} that the difference $\ord_{1}-\ord_{0}$ is uniformly bounded by elements of $\Gamma$: this property extends to the function $\ord_{1,\ult}-\ord_{0,\ult}$ on $R_{0,\ult}$. Since $\Gamma\subset C_{1}$, the last condition is therefore equivalent to $\ord_{1,\ult}\,(x_{\ult})>C_{1}$, hence to $x_{\ult}\in P_{C_{1}}$, which completes the proofs of \ref{domain} and \ref{regular}.\qed

\Subsection{Proof of the lifting theorem \ref{mainbis}}
The following proposition and its proof are essentially taken from \cite[Lemma 2.2]{BeDeLivdD79}.
\begin{subprop}\label{lift} 
Consider a commutative diagram of integral domains
$$\xymatrix{A_{\mathstrut}\ar@{}[r] |{\subset_{\mathstrut}}\ar@{^{(}->}[d]^{i}& V_{\mathstrut}\ar@{>>}^{\pi}[d]\\
A' \ar@{}[r] |{\subset} & V/P
}$$
where:
\begin{sitemize}
\item $V$ is a Henselian valuation ring, $P$ is a prime ideal of $V$ and $\pi$ is the canonical map;
\item $i$ is injective and the extension $\Frac(A')/\Frac(A)$ admits a separating transcendence basis.
\end{sitemize}
Then $A'$ lifts to $V$, i.e. there is a subring of $V$ containing $A$ and mapping isomorphically to $A'$ by $\pi$.
\end{subprop}

\pf put $F=\Frac(A)$ and $F'=\Frac(A')$. First, we may replace $A'$ by $F'\cap (V/P)$, which is a valuation ring because $V/P$ is. If $B$ is a separating transcendence basis for $F'/F$, then for each $b\in B$ we have $b\in A'$ or $b^{-1}\in A'$. So, by modifying $B$ we may assume that $B\subset A'$. Now the ring $A[B]$ lifts trivially to $V$ (just lift $B$ arbitrarily), so we assume from now on that $F'$ is separably algebraic over $F$. By Zorn's lemma, we are reduced to the case $A'=A[x]$ where $x$ is a root of $g\in A[X]$, irreducible and separable over $F$. So we have $g(x)=0$ and $g'(x)\neq0$. Let $\til{x}\in V$ be a lift of $x$: we have $g(\til{x})\in P$ and  $g'(\til{x})\notin P$, whence   $g'(\til{x})^2\notin P$ since $P$ is prime. So $e:=g(\til{x})/g'(\til{x})^2$ belongs to the maximal ideal of $V$. By the ``Hensel-Rychlik lemma'' (following from the Hensel property applied to the polynomial $G(h)=\frac{1}{g(\til{x})}\,g\left(\til{x}+e\,g'(\til{x})\,h\right)$) there exists $\ol{x}\in V$ with $g(\ol{x})=0$ and $\ol{x}\equiv\til{x}\bmod {e\,g'(\til{x})}$. In particular we have $\pi(\ol{x})=\pi(\til{x})=x$. Put $\ol{A}:=A[\ol{x}]\subset V$: then $\ol{A}$ lifts $A'$, because $\pi(\ol{A})=A[x]=A'$ and $A'$ and $\ol{A}$ can both be seen as subrings of $F[X]/(g(X))$.\qed

\begin{subcor}\label{liftpoints} With $A\subset V\surj V/P$ as in \rref{lift}, assume that the composite map $A\to V/P$ is injective and that the extension $\Frac(V/P)/\Frac(A)$ is \emph{separable}. Let $Y$ be an $A$-scheme locally of finite type. Then the natural map $Y(V)\to Y(V/P)$ is onto. 
\end{subcor}
\pf Since $V/P$ is a local ring, every morphism $y:\Spec\,(V/P)\to Y$ factors through an affine open subset of $Y$. So we may assume that $Y=\Spec\,(B)$ with $B$ finitely generated over $A$. Then $y$ corresponds to $\varphi:B\to V/P$. If $A'\subset V/P$ is the image of $\varphi$, then $\Frac(A')/\Frac(A)$ is a finitely generated separable extension, hence admits a separating transcendence basis \cite[V, \S\,9, $\mathrm{n}^\circ$\, 3, th.~2]{BourA4-5}. The conclusion then follows from  \ref{lift}. \qed

\begin{subrem}
Another noteworthy special case of \ref{lift} (already mentioned in \cite{BeDeLivdD79}) is when $\Frac(V/P)$ (hence also $\Frac(V)$) has characteristic zero: we may then take $A=\ZZ$ and $A'=V/P$ and conclude that $\pi$ has a section.
\end{subrem}
\Subsubsection{End of the proof of Theorem \ref{mainbis}. }
With $R$, $X$ and $C$ as in the theorem, we deduce from \ref{regular}\,\ref{regular2}  that  $\Frac\left(R^\ang{C}\right)$ is separable over $K$. Therefore we may apply Corollary \ref{liftpoints} with $A=R$, $V=R_{\ult}$, $P=P_{C}$ and $Y=X$. This completes the proof.

%%%%%%%
\section{Application: a closed image theorem}\label{sec:dempropre}
\Subsection{Basic topological facts}\label{TopFacts}
Recall that if $F$ is any Hausdorff topological field, we can uniquely define a topology on $X(F)$ for every $F$-scheme  $X$ locally of finite type, in such a way that ($X$ and $Y$ denoting arbitrary $F$-schemes locally of finite type):
\begin{sitemize}
\item if $X=\Aa^1_{F}$ we obtain the given topology on $X(F)=F$;
\item every $F$-morphism $f:X\to Y$ gives rise to a continuous map $X(F)\to Y(F)$ which, moreover, is an open (resp.\ closed)  topological embedding if $f$ is an open (resp.\ closed) immersion;
\item the natural bijection $(X\times Y)(F)\to X(F)\times Y(F)$ is a homeomorphism.
\end{sitemize}

In the sequel we keep the notations ($R$, $K$, $\Gamma$, $\ord$) of \ref{notations} and we take $F=K$ with the topology defined by the valuation. Thus, if $X$ is a $K$-scheme locally of finite type, we can characterize the topology on $X(K)$ as follows: for $x\in X(K)$, fix an affine open neighborhood $U=\Spec\,(A)$ of $x$ in $X$ and a finite sequence $(f_{1},\dots,f_{n})$ generating $A$ as a $K$-algebra. We obtain a basis of neighborhoods of $x$ in $X(K)$ by taking the ``balls'' $B(x,\gamma)=\left\{y\in U(K)\,\vert\, \ord\left(f_{i}(x)-f_{i}(y)\right)\geq\gamma, \:i=1,\dots, n\right\}$ for all $\gamma\in\Gamma$. 

Note that in the above description, $B(x,\gamma)$ is the image of $\cU(R)$ in $U(K)$, where $\cU$ is the spectrum of the $R$-algebra $\cA=R[\frac{f_{1}-f_{1}(y)}{t},\dots,\frac{f_{n}-f_{n}(y)}{t}]\subset A$ and we denote by $t$ any element of $K$ with valuation $\gamma$. (More generally, it can be checked that if $X$ is of finite type over $K$, we obtain a basis \modif{of open sets for  $X(K)$} by taking the sets $\im\left(\cX(R)\to X(K)\right)$ where $\cX$ runs through all $R$-schemes of finite type with generic fiber $X$).

If $\cX$ is a \emph{separated} $R$-scheme of finite type, then we can identify $\cX(R)$ with a subset of $\cX(K)=\cX_{K}(K)$, which is easily seen to be open; we can then endow  $\cX(R)$ with the induced topology. It is in fact possible to define the topology on  $\cX(R)$ directly, even if $\cX$ is not separated; however, this takes some more care (see \cite[Proposition 3.1]{ConAdelic}) and the present definition will be suffficient for our purposes. 

With $\cX$ as above, denote by $\cX_{0}$ the Zariski closure of $\cX_{K}$ in $\cX$, with its reduced subscheme structure. Then  $\cX_{0}$ and $\cX$ have the same $R$-points (resp.\ $K$-points), and it is easy to see that $\cX_{0}$ is flat over $R$ (recall that every torsion-free $R$-module is flat). It is also of finite type, as a closed subscheme of $\cX$, and hence of \emph{finite presentation} by \cite[(3.4.7)]{GruRa71}. To summarize, when using ``$R$-models'' to study the topology of a given $K$-scheme of finite type $X$, we only need models of $X_{\mathrm{red}}$ which are flat of finite presentation over $R$.\medskip

The following result is essentially equivalent to Corollary \ref{PHI} (in a more general setting, see also \cite[Proposition 4.1.1]{NN}):
\begin{prop}\label{EquivPHI} Assume that $R$ is Henselian and $\wh{K}$ is separable over $K$, and let $f:\cX\to\cY$ be a morphism of $R$-schemes of finite presentation, with $\cY$ separated. 

Then the induced map $f_{R}:\cX(R)\to \cY(R)$ has closed image.
\end{prop}
\dem The question is local on $\cY$, so we may assume that $\cY$ is affine, and even that $\cY=\Aa^n_{R}=\Spec\, R[T_{1},\dots, T_{n}]$ for some $n$, by choosing a closed immersion $\cY\inj\Aa^n$. Using a finite affine open covering of $\cX$, we may assume that 
$$\cX=\Spec \left(_{\mathstrut} R[T_{1},\dots, T_{n},Z_{1},\dots, Z_{m}]/(F_{1},\dots,F_{r})\right)$$
 for suitable polynomials $F_{j}\in R[T,Z]$. We may further assume that the origin $0\in \cY(R)=R^{n}$ is in the closure of the image of $f_{R}$. 

This means the following: for each $\gamma\in\Gamma$, there exist $t_{1},\dots, t_{n},z_{1},\dots, z_{m}$ in $R$ such that $F_{j}(t,z)=0$ ($1\leq j\leq r$) and $\ord\,(t_{i})\geq\gamma$ ($1\leq i\leq n$). Since $F_{j}$ has coefficients in $R$, this implies $\ord\,(F_{j}(0,z))\geq\gamma$. In other words, the fibre $\cY_{0}$ of $f$ at $0$ (which is an $R$-scheme of finite presentation) has $R_{\gamma}$-valued points for all $\gamma\in\Gamma^{+}$. By \ref{PHI}, \modif{$\cY_{0}(R)\neq\emptyset$}. In other words, $0$ is in the image of $f_{R}$.\qed 
\begin{subrem}
We have assumed $\cY$ separated only to avoid using the general definition of the topology on $\cY(R)$, alluded to in \ref{TopFacts} above. Surprisingly (at least to the author), this assumption is not necessary, and in fact $\cY(R)$ is always a Hausdorff space, even if  $\cY$ is not separated.
\end{subrem}
\Subsection{Proof of Theorem \ref{ThPropre}}
Consider $f:X\to Y$ as in Theorem \ref{ThPropre}. To prove that the image of $f_{K}$ is closed, we may assume $Y$ affine. Fix a flat, affine, finitely presented $R$-scheme $\cY$ with generic fiber $Y$. It suffices to show that $f_{K}(X(K))\,\cap\, \cY(R)$ is closed in $\cY(R)$, because the sets $\cY(R)\subset Y(K)$, for varying $\cY$, form a basis of open subsets. 

Since $K$ is the increasing union of its subrings $R[t^{-1}]$, where $t$ runs through nonzero elements of $R$, we can apply the results of \cite[\S 8]{EGA4_III} and find, for suitable such $t$, a scheme $\cX_{1}$, separated  of finite presentation over $R[t^{-1}]$ (hence also over $R$), such that $(\cX_{1})_{K}=X$, and an $R$-morphism $f_{1}:\cX_{1}\to \cY$ extending $f$. 

By Nagata's compactification theorem \cite[Theorem 4.1]{ConNagata} the morphism $f_{1}$ factors as  $\cX_{1}\mathop{\hookrightarrow}\limits^{j}\cX\xrightarrow{\ol{f}}\cY$ where $j$ is a dense open immersion and $\ol{f}$ is proper; in particular, we have $\cX_{K}=\left(\cX_{1}\right)_{K}=X$ since $\left(\cX_{1}\right)_{K}$ is \modif{assumed} proper over $\cY_{K}$. Thus, by the valuative criterion of properness, we have  $f_{K}(X(K))\,\cap\, \cY(R)=\ol{f}(\cX(R))$. Hence the result follows from \ref{EquivPHI}, applied to $\ol{f}$.\qed
\begin{subrem}
Of course, if $f$ is assumed projective, Nagata's theorem is not needed: in the proof, we can choose $\cX_{1}$ quasiprojective over $R$. 
\end{subrem}

\bibliography{GVbiblio}

 \end{document}